\documentclass[10pt,reqno]{amsart}
\usepackage{amsmath, latexsym, amsfonts, amssymb,amsthm, amscd,epsfig,enumerate}
\usepackage{subfigure,color, float, caption}
\usepackage{tikz,colortbl}
\usetikzlibrary{calc}
\usetikzlibrary{decorations.pathreplacing}

\advance\hoffset -.75cm

\definecolor{refkey}{gray}{.75}
\definecolor{labelkey}{gray}{.75}

\newcommand{\R}{\mathbb R}
\newcommand{\Z}{\mathbb Z}
\newcommand{\N}{\mathbb N}

\newcommand{\pr}{\mathbb P}

\newtheorem{teo}{Theorem}[section]
\newtheorem{Theorem}{Theorem}[section]
\newtheorem{Lemma}[teo]{Lemma}
\newtheorem{Corollary}[teo]{Corollary}
\newtheorem{Remark}[teo]{Remark}
\newtheorem{Proposition}[teo]{Proposition}
\newtheorem{Definition}[teo]{Definition}
\newtheorem{Example}[teo]{Example}

\newtheorem{Assumption}[teo]{Assumption}

\title
{On the critical parameters of branching random walks}

\author[D.~Bertacchi]{Daniela Bertacchi}
\address{D.~Bertacchi, Dipartimento di Matematica e Applicazioni,
Universit\`a di Milano--Bicocca,
via Cozzi 53, 20125 Milano, Italy.}
\email{daniela.bertacchi\@@unimib.it}

\author[F.~Zucca]{Fabio Zucca}
\address{F.~Zucca, Dipartimento di Matematica,
Politecnico di Milano,
Piazza Leonardo da Vinci 32, 20133 Milano, Italy.}
\email{fabio.zucca\@@polimi.it}

\date{}

\begin{document}

\begin{abstract}Given a discrete spatial structure $X$, we define continuous-time branching processes
	$\{\eta_t\}_{t \ge 0}$
	that model a population breeding and dying on $X$. These processes are usually called branching random walks,
	and $\eta_t(x)$ denotes the number of individuals alive at site $x$ at time $t$.
	They are characterized by breeding rates $k_{xy}$ (governing the rate at which individuals at $x$ send offspring to $y$), and by a multiplicative speed parameter $\lambda$. 
	These processes also serve as models for epidemic spreading, where $\lambda k_{xy}$ represents the infection rate from $x$ to $y$.
	In this context, $\eta_t(x)$ represents the number of infected individuals at $x$ at time $t$, and the removal of an individual is due to either death or recovery.
	Two critical parameters of interest are the global critical parameter $\lambda_w$, related to global survival, and the local critical parameter $\lambda_s$, related to survival within finite sets (with $\lambda_w\le\lambda_s$).
	In disease or pest control, the primary goal is to lower $\lambda$ so that the process dies out, at least locally. Nevertheless, a process that survives globally can still pose a threat, especially if sudden changes cause global survival to transition into local survival. In fact, local
	modifications to the rates can affect the values of both critical parameters,  making it important to understand when and how they can be increased.
	Using results on the comparison of extinction probabilities for a single branching random walk across different sets, we extend the analysis to extinction probabilities and critical parameters of pairs of branching random walks whose rates coincide outside a fixed set $A \subseteq X$.
	We say that two branching random walks are equivalent if their rates coincide everywhere except on a finite subset of $X$. Given an equivalence class of branching random walks, we prove that if one process has $\lambda^*_w \neq \lambda^*_s$, then $\lambda^*_w$ is the maximal possible value of this parameter within the class.
	We describe the possible configurations for the critical parameters within these equivalence classes.
\end{abstract}

\maketitle
\noindent {\bf Keywords}: branching random walk, branching process, critical parameters, local survival, global survival, pure global survival phase.

\noindent {\bf AMS subject classification}: 60J05, 60J80.

\section{Introduction}

A branching process, or the Galton-Watson process (see \cite{cf:GW1875}), is a process where a particle dies and gives birth to a random number of offspring,
according to a given offspring law $\rho$ ($\rho(n)$ being the probability of having exactly $n$ children).
Different particles breed independently, according to $\rho$. 
The process can either go extinct (i.e. no particles are alive from a certain time on) almost surely,
or survive (i.e. at least one particle in any generation) with positive probability.
The probability of extinction can be computed as a function of $\rho$.

A branching random walk (BRW hereafter) is a process where particles are described by their location $x\in X$,
where $X$ is an at most countable set ($X$ is usually interpreted as a spatial variable, but can also be seen
as a ``type'', see for instance \cite{cf:KurtzLyons}). The particles at site $x\in X$ are replaced by a random
number of children, which are placed at various locations on $X$. The breeding law depends on the site where the parent lives, and
all particles breed independently.
This class of processes (in continuous and
discrete time) has been studied by many authors: see \cite{cf:AthNey, cf:Big1977, cf:Harris63} for older results.
Besides the classical environment $\Z^d$, the process has been studied in 
other settings such as finite sets
\cite{cf:MountSchi}, or trees \cite{cf:Attia, cf:Ligg1, cf:MadrasSchi, cf:Muller1, cf:PemStac1, cf:Su}.
Note that, in the case of the branching random walk,
there is no upper bound for the number of particles per site. When such an upper bound is fixed, say at most
$m$ particles per site, we get the $m$-type contact process. The branching random walk can be obtained
as the limiting process as $m$ goes to infinity (\cite{cf:PemStac1,cf:Z1}).

BRWs can be considered in discrete time, where each particle is alive for one generation and its offspring all live in the next generation, or
in continuous time, where generations overlap since particles may breed at different times during their whole lifetime.
In particular, one can consider a family of parametrized processes, where we fix the breeding rates between locations and the parameter $\lambda>0$
represents the reproductive speed (the larger $\lambda$, the shorter the time intervals between subsequent breedings, see Section \ref{sec:basic} for details).

Since the BRW lives on a spatial structure, its behaviour is in general more complex than the one of a branching process: 
indeed survival and extinction can be studied globally (on the whole space $X$) but also locally (on a single site or a finite set $A$).
Starting the process with one particle at a given site $x$, and fixing a set $A\subseteq X$ only one of the following holds for the BRW:
(1) it goes almost surely extinct, (2) it survives globally but not locally,
(3) it survives globally and locally but with different probabilities (non-strong local survival),
(4) it survives globally and locally with equal probability (strong local survival).
We stress that no strong local survival means that there is either non-strong local survival or almost sure local extinction.

When considering BRWs in continuous time, a fixed family may exhibit different behaviours for different values of the speed parameter $\lambda$,
and the threshold that separates global extinction from global survival is called global critical parameter $\lambda_w$.
Analogously, the threshold between local extinction and local survival is called local critical parameter $\lambda_s$.
We note that the subscripts $w$ and $s$ stand for "weak" and "strong" which are sometimes used as synonyms for "global"
and "local", respectively.

Many authors have addressed natural questions, such as how to identify the critical parameters (\cite{cf:PemStac1,cf:Stacey}),
or establish criteria for survival and extinction (\cite{cf:Gantert1, cf:Gantert2, cf:MachadoMenshikovPopov, cf:Muller1, cf:Muller2, cf:Su}).
Once survival with positive probability is established, one can address the problem of identifying and computing extinction probabilities, see
\cite{cf:Haupt}.
In ecological or epidemiological contexts, the goal may be either to promote survival (e.g., of an endangered species) 
or to achieve extinction (e.g., for disease or pest control). 
In the case of continuous-time branching random walks, two key factors are involved: the speed parameter $\lambda$ and 
the family of rates $\{k_{xy}\}_{x,y \in X}$. 
If the goal is to drive an epidemic to extinction, one strategy is to lower $\lambda$ (for example, wearing a mask to reduce respiratory infections 
such as COVID-19 serves this purpose). When $\lambda \le \lambda_s$, we know that the infection will eventually disappear from any site. 
However, if $\lambda_w < \lambda$, it may still persist in the system.
Another strategy is to reduce the infection rates $k_{xy}$, as doing so may increase the values of the critical parameters and 
thereby change the behaviour of the process for a fixed $\lambda$.
More precisely, suppose that $\lambda$ is fixed and that $\lambda_w$, $\lambda_s$ are the critical parameters associated with the rates $\{k_{xy}\}_{x,y \in X}$.
After some modifications, the rates become $\{k^*_{xy}\}_{x,y \in X}$, possibly leading to different critical parameters $\lambda_w^*$ and $\lambda_s^*$.
If $\lambda_s <\lambda <\lambda_s^*$, then the modified process dies out in all locations, even though the original process did not.
Altering the infection rates corresponds to changing the interactions between sites, for example, by implementing social distancing or closing schools. 
However, in practice, such changes are typically feasible only on finite subsets of $X$ (i.e., they are local changes).
It is therefore of utmost importance to understand how local changes affect the critical parameters, as our goal is to actively influence the course of an epidemic through these interventions.
Another important consideration arises when the process is in a purely global survival phase ($\lambda_w <\lambda \le\lambda_s$): the epidemic appears "mild," since it disappears from all locations. However, local modifications can alter this balance, and if they reduce the local critical parameter (i.e., if $\lambda_s^* < \lambda$), the disease could reappear at all sites.
Our results show that it is impossible to increase $\lambda_w$ when the process is in a pure global survival phase. However, in other cases, it is possible to do so. For instance, isolating superspreaders may be effective, as illustrated in Example~\ref{exmp:modifiedtree}.
In that example, there is a single superspreader location, the origin of the tree, denoted by $o$. When $k_{oo}$ is large and $\lambda \in (0, 1/d)$, the process exhibits local survival. 
However, if $k^*_{oo}$ is sufficiently small (i.e., the superspreader is isolated), the process becomes extinct both locally and globally.

The present paper is devoted to the question of what can be said when a BRW is modified locally, that is, when the breeding laws are changed only within a fixed set $A \subseteq X$.
For instance, one may ask how the modified BRW behaves if the original BRW does not exhibit strong local survival. 
Theorem~\ref{th:modifiedBRW} shows that there is global survival and no strong local survival in $A$ for the original BRW if and only if the same holds for the modified BRW, regardless of the specific modifications made in $A$.
As a corollary, we obtain that if the original BRW dies out locally in $A$ and the modified BRW survives globally, then almost sure global extinction of the original BRW is equivalent to strong local survival in $A$ for the modified one (see Corollary~\ref{cor:pureweak-nonstrong2} and Figure~\ref{fig:modifiedBRW}).
Moreover, for a fixed irreducible BRW, if there is global survival and no strong local survival in some set $A \subseteq X$, then this property holds for every finite set $B \subseteq X$.

Building on these results in discrete time, we show that in continuous time, a local modification of the BRW in a finite subset $A$, which lowers the global critical parameter (typically achieved by adding a sufficiently fast reproduction rate at some site), implies that the global and local critical parameters of the modified BRW coincide (see Corollary~\ref{cor:pureweak-nonstrong}).
This leads to a general method for constructing examples such as Example~\ref{exmp:modifiedtree}, where the modified BRW exhibits strong local survival for some values of the parameter below a certain threshold and above another, and non-strong local survival for intermediate values (see Figures~\ref{fig:modifiedBRW} and \ref{fig:modifiedtree}).
This example was originally presented in \cite{cf:BZ14-SLS}, but it is included here with a simpler proof and within a more general framework.
Furthermore, we prove that in general, a continuous-time BRW obtained via a local modification that lowers the global critical parameter dies out globally at its critical value, a property that does not always hold (see \cite[Example 3]{cf:BZ2}).

Here is an outline of the paper. In Section~\ref{sec:basic}, we recall the definitions of discrete-time and continuous-time BRWs, the extinction probability vectors, and their asymptotic behaviour, specifically, whether extinction occurs as time tends to infinity. We also define global, local, and strong local survival, along with the critical parameters associated with continuous-time BRWs.
Theorem~\ref{th:spatial} and Corollary~\ref{cor:sen2} present known results about approximating a BRW by a sequence of BRWs, stated in a form that will be useful later. Proposition~\ref{pro:lambdaequal} provides a sufficient condition for the equality $\lambda_w = \lambda_s$.
In Section~\ref{sec:survivalprob}, we first compare extinction probabilities of a single BRW restricted to different sets (Theorem~\ref{th:strongconditioned}), and then compare extinction probabilities of two different BRWs (Theorem~\ref{th:modifiedBRW}).
Next, we study how the critical parameters of a continuous-time BRW are affected by local modifications of its rates (see Corollaries~\ref{cor:pureweak-nonstrong}, \ref{cor:pureweak-nonstrong2}, and \ref{cor:pureweak-nonstrong3}).
Section~\ref{sec:max} begins with the observation that if a BRW has a pure global survival phase ($\lambda_w < \lambda_s$, meaning there exists a range of values of the speed parameter $\lambda$ where the process survives globally but eventually leaves all finite sets), then no local modification of the rates can increase the global critical parameter (see Proposition~\ref{pro:maximality}).
This implies that the global critical parameter of a BRW with a pure global survival phase is maximal among those of all BRWs obtained through local modifications.
%
%
%
Finally, in Section~\ref{sec:max}, we analyse the behaviour of the BRW on a homogeneous tree with one added reproductive loop at the origin (Example~\ref{exmp:modifiedtree}), and address several general questions about the critical parameters when the rates are modified in a subset $A$.

\section{Basic definitions and preliminaries}
\label{sec:basic}

\subsection{Discrete-time Branching Random Walks and their survival/extinction.}
Given an at most countable set $X$, we define a discrete-time BRW
as a process $\{\eta_n\}_{n \in \N}$,
where $\eta_n(x)$ is the number of particles alive at $x \in X$ at time $n$. 
The dynamics is described as follows: 
let $S_X:=\{f:X \to \N\colon \sum_yf(y)<\infty\}$ and let 
$\mu=\{\mu_x\}_{x \in X}$ be a family of probability measures
on the (countable) measurable space $(S_X,2^{S_X})$. 
A particle of generation $n$ at site $x\in X$ lives one unit of time;
after that, a function $f \in S_X$ is chosen at random according to the law $\mu_x$.
This function describes the number of children and their positions, that is,
the original particle is replaced by $f(y)$ particles at
$y$, for all $y \in X$. The choice of $f$ is independent for all breeding particles.
The BRW is denoted by $(X,\mu)$.

Some results rely on the \textit{first moment matrix}
$M=(m_{xy})_{x,y \in X}$,
where each entry
$m_{xy}:=\sum_{f\in S_X} f(y)\mu_x(f)$ represents
the expected number of children that a particle living
at $x$ sends to $y$
(briefly, the expected number of particles from $x$ to $y$).
For the sake of simplicity, we require
that $\sup_{x \in X} \sum_{y \in X} m_{xy}<+\infty$.

To  a generic discrete-time BRW we associate a graph $(X,E_\mu)$, where $(x,y) \in E_\mu$  
if and only if $m_{xy}>0$.
In contrast, given a graph $(X,E)$, we say that a BRW $(X,\mu)$ is adapted to the graph if 
$m_{xy}>0$ if and only if $(x,y)\in E$.
We say that there is a path of lenght $n$ from $x$ to $y$, and we write $x \stackrel{n}{\to} y$, if it is
possible to find a finite sequence $\{x_i\}_{i=0}^n$ (where $n \in \N$)
such that $x_0=x$, $x_n=y$ and $(x_i,x_{i+1}) \in E_\mu$
for all $i=0, \ldots, n-1$ (observe that there is always a path of length $0$ from $x$ to itself). If $x \stackrel{n}{\to} y$ for some $n \in \mathbb{N}$,
then we write $x \to y$; whenever $x \to y$ and $y \to x$ we write $x \rightleftharpoons y$.
The equivalence relation $\rightleftharpoons$ induces a partition of $X$: the
class $[x]$ of $x$ is called \textit{irreducible class of $x$}.
If the graph $(X,E_\mu)$ is \textit{connected} (that is, there is only one irreducible class),
then we say that the BRW is \textit{irreducible}.
Irreducibility implies that the progeny of any particle can spread to any site of the graph with positive probability.

We consider initial configurations with only one particle placed at a fixed site $x$ and we
denote by $\pr^{x}$ 
the law of the corresponding process. Evolution of processes with more than one initial particle
can be obtained by superimposition.

In the following, \textit{wpp}
is shorthand for ``with positive probability'' (although, when talking about survival, 
``wpp'' will usually be tacitly understood).
In order to avoid trivial situations where particles have one offspring almost surely, we assume
henceforth the following.
\begin{Assumption}\label{assump:1}
	For all $x \in X$ there is a vertex $y \rightleftharpoons x$ such that
	$\mu_y(f\colon  \sum_{w\colon w \rightleftharpoons y} f(w)=1)<1$,
	that is, in every equivalence class (with respect to $\rightleftharpoons$)
	there is at least one vertex where a particle
	can have, inside the class, a number of children different from 1 wpp.
\end{Assumption}
We now introduce some definitions. The notation of the main quantities is summarized in Table \ref{tb:table1}.
For a more formal statement, see Definition \ref{def:extprob} and Equation \eqref{eq:criticalparameters}.
\begin{table}
	\begin{tabular}{ |c|c| }
		\hline 
		\rowcolor{lightgray}
		Notation & Meaning  \\ \hline
		${\mathbf{q}}(x,A)$ &  probability of extinction in $A$, starting from $x$\\ 
		\hline
		${\mathbf{q}}(A)$ &  vector whose components are ${\mathbf{q}}(\cdot,A)$ \\ 
		\hline
		$\lambda_w(x)$ & threshold between global extinction and global survival, starting from $x$ \\ 
		\hline
		$\lambda_s(x)$ & threshold between extinction and survival in $x$, starting from $x$   \\ 
		\hline
		\end{tabular}
	\captionof{table}{Main notation. The BRW starts with one individual in $x\in X$.}\label{tb:table1}
\end{table}

\begin{Definition}\label{def:extprob}
	Given a BRW $(X,\mu)$, $x\in X$ and $A\subseteq X$, the probability of extinction in $A$, starting from a particle at $x$, is defined as
	\[
	{\mathbf{q}}(x,A)
	:=1-\pr^{x}
	(\limsup_{n \to \infty} \sum_{y \in A} \eta_n(y)>0).
	\]
	We denote by ${\mathbf{q}}(A)$ the extinction probability vector, whose $x$-entry is ${\mathbf{q}}(x,A)$.
\end{Definition}
If $A=\{y\}$, we write ${\mathbf{q}}(x,y)$ instead of ${\mathbf{q}}(x, \{y\})$.
Note that
${\mathbf{q}}(x,A)$ depends on $\mu$. When we need to stress this dependence, we write 
$ {\mathbf{q}}^\mu(x,A)$.
Extinction probabilities have been the object of intense study during the last decades.
We refer the reader, for instance, to \cite{cf:BBHZ, cf:BZ4, cf:Z1}.

\begin{Definition}\label{def:survival} 
	Let $(X,\mu)$ be a BRW starting from one particle at $x\in X$ and let $A\subseteq X$.
	We say that
	\begin{enumerate}
		\item 
		the process \textsl{survives locally wpp} in $A$ 
		if 
		$
		{\mathbf{q}}(x,A)<1
		$;
		\item
		the process \textsl{survives globally wpp} if
		$
		{\mathbf{q}}(x,X) <1$;
		\item
		there is \textsl{strong local survival wpp} in $A$ 
		if
		$ 
		{\mathbf{q}}(x,A)={\mathbf{q}}(x,X)<1$
		and \textsl{non-strong local survival wpp} in $A$ if ${\mathbf{q}}(x,X)<{\mathbf{q}}(x,A)<1$;
		\item 
		there is \textsl{no strong local survival} in $A$
		if either ${\mathbf{q}}(x,A)=1$ or ${\mathbf{q}}(x,X)<{\mathbf{q}}(x,A)$;
		\item
		the process is in a \textsl{pure global survival phase} if
		$
		{\mathbf{q}}(x,X) <{\mathbf{q}}(x,x)=1
		$.
	\end{enumerate}
\end{Definition}
When there is no survival wpp, we say that there is extinction
and the fact that extinction occurs 
almost surely, will be tacitly understood.
When there is strong local survival, it means that for almost all realizations the process either survives locally
(hence globally) or it goes globally extinct. More precisely,
there is strong survival at $y$ starting from $x$ if and only if the probability
of local survival at $y$ starting from $x$ conditioned on global survival starting from $x$ is $1$.

We want to stress that, for a BRW starting from $x\in X$,
${\mathbf{q}}(x,X)={\mathbf{q}}(x,A)$ if and only if global survival 
is
equivalent to strong local survival at $A$.
On the other hand,
${\mathbf{q}}(x,X)<{\mathbf{q}}(x,A)$
if and only if there is global survival and no strong local survival at $A$.

In general, Definition \ref{def:survival} depends on the starting vertex.
However, if the process is irreducible, then the process survives locally or
globally starting from $x$ if and only if the same holds when starting from $y$,
for all $x,y\in X$. Strong local survival may still depend on the starting vertex, 
even in the irreducible case. If at all sites we have a positive probability of
no children, then one can prove that strong local survival does not depend on the starting vertex (see \cite[Section 3]{cf:BZ14-SLS}).

Moreover, if a BRW is irreducible, $A,B \subset X$ are nonempty, finite subsets and $C \subset X$ is nonempty, then $\mathbf{q}(x,A)=\mathbf{q}(x,B) \ge \mathbf{q}(x,C)$ for all $x\in X$. The first equality holds since local survival does not depend on the target vertex and survival in a finite (nonempty) set is equivalent to local survival to a vertex in the set; more precisely, in the irreducible case, if a vertex $y$ is visited infinitely often, then the conditional probability of visiting infinitely many times all other vertices is 1 (by Borel-Cantelli Lemma).  The second inequality follows from the fact that any nonempty set contains a finite nonempty set.

Henceforth we make use of the natural partial order between vectors defined as follows: $\mathbf{q}(A) \le \mathbf{q}(B)$ if and only if $\mathbf{q}(x,A) \le \mathbf{q}(x,B)$ for all $x \in X$. Therefore, $\mathbf{q}(A) < \mathbf{q}(B)$ if and only if 
$\mathbf{q}(x,A) \le \mathbf{q}(x,B)$ for all $x \in X$ and $\mathbf{q}(x_0,A) < \mathbf{q}(x_0,B)$ for some $x_0 \in X$;
moreover $\mathbf{q}(A) \nleq \mathbf{q}(B)$ if and only if $\mathbf{q}(x_0,A) > \mathbf{q}(x_0,B)$ for some $x_0 \in X$.

\subsection{Continuous-time Branching Random Walks}
\label{subsec:continuous}

Given an at most countable $X$ and a nonnegative matrix $K=(k_{xy})_{(x,y)\in X\times X}$, 
one can define a family of continuous-time Branching Random Walks
$\{\eta_t\}_{t\ge0}$, where $\eta_t(x)$ represents the number of particles alive at time $t$ at site $x$, for any $x\in X$.
The family is indexed by the reproductive speed parameter $\lambda>0$.

Each particle has an exponentially distributed
lifetime with parameter 1.
During its lifetime each particle alive at $x$
breeds into $y$ according to the arrival times of its own Poisson process with
parameter $\lambda k_{xy}$ (representing the reproduction rate).
We denote by  $(X,K)$ this family of continuous-time BRWs (depending on $\lambda>0$).
It is not difficult to see that the introduction
of a nonconstant death rate $\{d(x)\}_{x \in X}$ does not represent a
significant generalization. In fact,
one can study 
a new BRW with death rate 1 and
reproduction rates $\{\lambda k_{xy}/d(x)\}_{x,y \in X}$; the two processes have the same behaviours in
terms of survival and extinction (\cite[Remark 2.1]{cf:BZ14-SLS}).

To each continuous-time BRWs 
we associate a discrete-time counterpart, namely $(X,\mu)$ where we simply take into account the number and the positions of all the offsprings
born before the death of the parent.
Clearly, from the knowledge of the discrete-time counterpart of $(X,K)$, we cannot retrieve how many particles are alive at time $t$ (and where),
but the information about survival and extinction is intact.
The probabilities of extinction of a continuous BRW coincides with those of its discrete-time counterpart and they depend on $\lambda$: we denote them by $\mathbf{q}(x,A|\lambda)$ although the dependence on $\lambda$ will be usually omitted.   
In particular, we extend every definition from the discrete-time case to the continuous-time case in a natural way by using the discrete-time counterpart of a continuous-time process. More precisely, when we say that a continuous-time process has a certain property, we mean that its discrete-time counterpart has it; for instance, a continuous-time process is irreducible (by definition) if and only if its discrete-time counterpart is.
We note that, in the discrete-time counterpart of a 
continuous-time BRW, 
at every vertex there is a positive probability of dying without  breeding; hence, in the continuous-time irreducible case, either strong local survival is a common property of all starting vertices, or it holds nowhere.
Moreover, it is easy to show that the expected number of children from $x$ to $y$, for the discrete-time counterpart of the process, is $m_{xy}=\lambda k_{xy}$.

Given $x \in X$, two critical parameters are associated with the
continuous-time BRW: the \textit{global} 
\textit{survival critical parameter} $\lambda_w(x)$ and the  \textit{local} 
\textit{survival critical parameter} $\lambda_s(x)$ (or, briefly, the global and local critical parameters, respectively) defined as
\begin{equation}\label{eq:criticalparameters}
	\begin{split}
		\lambda_w(x)&:=\inf \Big \{\lambda>0\colon \,
		\pr^{x}
		\Big (\sum_{w \in X} \eta_t(w)>0, \forall t\ge 0\Big) >0 \Big \},\\
		\lambda_s(x)&:=
		\inf\{\lambda>0\colon \,
		\pr^{x} \big(\limsup_{t \to \infty} \eta_t(x)>0 \big) >0
		\}.
	\end{split}
\end{equation}
These values depend only on the irreducible class of $x$. In particular, they are constant
if the BRW is irreducible, 
in which case we simply write $\lambda_w$ and $\lambda_s$.
Indeed, suppose that $x \to y$; since there is a positive probability that a descendant of a particle in $x$ is placed in $y$, if there is a positive probability of survival in a set $A$ starting from $y$, then the same holds starting from $x$. Therefore, taking $A=X$, $\lambda_w(x) \le \lambda_w(y)$ and if $x \rightleftharpoons y$, then $\lambda_w(x) = \lambda_w(y)$. Moreover, if $x$ visited an infinite number of times then, by Borel-Cantelli's Lemma, an infinite number of descendants is placed at $y$: more precisely, the probability of visiting $y$ infinitely often conditioned on visiting $x$ infinitely often is 1. Thus, if $x \rightleftharpoons y$ then $\lambda_s(x) = \lambda_s(y)$. This last equality can be proven also by using equation~\eqref{eq:lambdas}; this simple exercise is left to the reader.

Note that, while in the discrete-time case there is only one process that can be in a pure global survival phase or not, 
in continuous-time we have a family of processes, indexed by $\lambda$. Hence we may have processes, in the same family,
in pure global survival, and others which are not in pure global survival.
We say that there exists a \textit{pure global survival phase} starting from $x$,
if the interval $(\lambda_w(x),\lambda_s(x))$ is nonempty.
No reasonable definition of a \textit{strong local survival critical parameter} is possible
(see \cite{cf:BZ14-SLS}).

It is possible to identify $\lambda_s(x)$ in terms of the  matrix $K$.
Indeed, consider the
$n$th power matrix $K^n$ with entries $k^{(n)}_{xy}$: it is possible to prove that
\begin{equation}\label{eq:lambdas}
	\lambda_s(x)=
	1/\limsup_{n \in \N} \sqrt[n]{k^{(n)}_{xx}} 
\end{equation}
(see \cite[Theorems 4.1 and 4.7]{cf:BZ2} or \cite[Lemma 3.1]{cf:PemStac1}). One can also write the local critical parameter in terms of a generating function.
Let  
\[\begin{split}
	\varphi^{(n)}_{xy}&:=\sum_{x_1,\ldots,x_{n-1} \in X \setminus\{y\}} k_{x x_1} k_{x_1 x_2} \cdots k_{x_{n-1} y}\\
	\Phi(x,y|\lambda)&:=\sum_{n =1}^\infty \varphi_{xy}^{(n)} \lambda^n.
\end{split}
\]
The coefficients $\varphi^{(n)}_{xy}$ play in the BRW theory
the same role played by the taboo probabilities in random walk theory: when $x=y$ these are called first-return probabilities
(see \cite[Section 1.C]{cf:Woess09}). Roughly speaking, $\lambda^n \varphi^{(n)}_{xy}$ is
the expected number of particles alive at $y$ at time $n$,
when the initial state is just one particle at $x$ and the
process behaves like a BRW except that every particle reaching
$y$ at any time $i <n$ is immediately killed (before breeding).
The following  characterization holds (see~\cite{cf:BZ2}):
\begin{equation}
	\label{eq:lambdas1}
	\lambda_s(x)=
	\max\{ \lambda 
	\in {\mathbb R}:\Phi(x,x|\lambda)\leq 1\}=
	\sup 
	\{ \lambda 
	\in {\mathbb R}:\Phi(x,x|\lambda)< 1\}
\end{equation}
where the last equality is due to the fact that $\lambda \mapsto \Phi(x,x|\lambda)$ is a continuous, strictly increasing power series from $[0, r)$ to $[0,+\infty)$ where $r$ is the radius of convergence (it is also left continuous in $r$ even if $\Phi(x,x|r)$ is infinite).
The characterization of $\lambda_w(x)$ is less explicit:
if we denote by $l_+^\infty(X)$ the usual subset of nonnegative, bounded functions on $X$,  we have a characterization in terms of the existence of solutions of certain inequalities
(see \cite[Theorem 4.2]{cf:BZ2})
\begin{equation}
	\label{eq:lambdaw}
	\lambda_w(x)
	=\inf \{\lambda \in \R: \exists v \in l^\infty_+(X), \,  v(x)> 0 \textrm{ and } \exists n \in \mathbb{N}\setminus \{1\},\, \lambda^n K^n v \ge v \}
\end{equation}
where $K^n v(x):=\sum_{y \in X} k^{(n)}_{xy}v(y)$ ($k^{(n)}_{xy}$ are the entries of the matrix $K^n$, the $n$th power of the matrix $N$).
Alternatively, there is global survival for a fixed $\lambda>0$ starting from $x$ if and only if there exists $v \in l^\infty_+(X)$ such that $v(x)>0$ and
$\frac{\lambda Kv(y)}{1+\lambda Kv(y)} \ge v(y)$ for all $y \in X$ (see \cite[Theorem 4.2]{cf:BZ2}).

The relations between the critical parameters and the extinction probabilities are summarized in Table~\ref{tb:table3}.

\begin{table}
	\begin{tabular}{ |c|c| }
		\hline 
		\rowcolor{lightgray}
		Relations between critical parameters and probabilities of extinction  \\ \hline
		$\lambda \le \lambda_s(x) \Longleftrightarrow {\mathbf{q}}(x,\{x\}|\lambda)= 1$ \\ 
		\hline
		$\lambda > \lambda_s(x) \Longleftrightarrow {\mathbf{q}}(x,\{x\}|\lambda)< 1$ \\ 
		\hline
			$\lambda < \lambda_w(x) \Longrightarrow {\mathbf{q}}(x,X|\lambda)= 1$, \quad ${\mathbf{q}}(x,X|\lambda)= 1 \Longrightarrow \lambda \le \lambda_w(x)$  \\ 
		\hline
			$\lambda > \lambda_w(x) \Longrightarrow {\mathbf{q}}(x,X|\lambda)< 1$, \quad ${\mathbf{q}}(x,X|\lambda)< 1 \Longrightarrow \lambda \ge \lambda_w(x)$  \\
		\hline
	\end{tabular}
	\captionof{table}{Relations between ctitical parameters and extinction probabilities; see \cite[Theorem~4.7 and Example~3]{cf:BZ2}.}\label{tb:table3}
\end{table}

\begin{Remark}
	\label{rem:finiteX}
	Even though we do not assume that the set $X$ is infinite, the finite case is straightforward. Indeed, if $X$ is finite and the BRW is irreducible, then $\lambda_s = \lambda_w$ coincides with the spectral radius of the nonnegative, irreducible matrix $K$ (by the Gelfand formula). It is easy to construct a solution to equation~\eqref{eq:lambdaw} for any $\lambda$ larger than the spectral radius of $K$: one can simply take a positive eigenvector corresponding to the Perron-Frobenius eigenvalue.  
	Thus, in the finite case, the critical values are given by the inverse of the Perron-Frobenius eigenvalue of $K$. From a probabilistic point of view, the equality $\lambda_s = \lambda_w$ follows from the fact that survival in a finite set is equivalent to visiting at least one site infinitely often; by irreducibility, this implies survival at every site.  
	
	More formally, the equality $\lambda_s = \lambda_w$ when $X$ is finite follows from \cite[Remark 4.1, Theorems 4.1 and 4.3]{cf:BZ2} (since $M_w^- = M_s$ in the finite irreducible case; see that paper for the necessary definitions) and from \cite[Theorem 3.2, Corollary 3.1, and the subsequent discussion]{cf:BZ14-SLS}.
\end{Remark}

\subsection{Spatial Approximation}

Given a sequence of BRWs $\{(X_n,\mu_n)\}_{n \in \N}$ such that $\liminf_{n \to \infty} X_n=X$,
we define $m(n)_{xy}:=\sum_{f \in S_{X_n}} f(y) \mu_{n,x}(f)$, where $x,y \in X_n$, and
the corresponding sequence of first moment matrices $\{M_n\}_{n \in \N}$. 
Clearly, for all $x, y \in X$, $m(n)_{xy}$ is eventually well defined as $n \to +\infty$. Indeed, $\liminf_{n \to \infty} X_n=X$, implies that given $A \subseteq X$ finite, $A \subseteq X_n$ for every sufficiently large $n$.
We wonder when a BRW can be approximated by a sequence of BRWs: more precisely, we are focusing on the local behaviour.
Note that in the following results we are not assuming that the BRW is irreducible. 

\begin{Theorem}[{\cite[Theorem 5.2]{cf:Z1}}]\label{th:spatial}
	Let us 
	fix a vertex $x_0 \in X$. Suppose that $\liminf_{n \to \infty} X_n=X$ and 
	$m(n)_{xy} \le m_{xy}$ for all $x,y \in X_n$ and all $n \in \N$. Assume that, for all $x,y \in X$, 
	$m(n)_{xy} \to m_{xy}$ as $n \to \infty$.
	\begin{enumerate}
		\item if
		$(X,\mu)$ dies out locally (resp.~globally) a.s.~starting from $x_0$, then $(X_n,\mu_n)$ dies out locally (resp.~globally)
		a.s~starting from $x_0$  for all $n \in \N$;
		\item if
		$(X,\mu)$ survives locally starting from $x_0$,  then $(X_n,\mu_n)$ survives locally
		starting from $x_0$ eventually as $n \to\infty$.
	\end{enumerate}
\end{Theorem}

Theorem~\ref{th:spatial} yields an analogous result for BRWs in continuous time.

\begin{Corollary}[{\cite[Corollary 5.3]{cf:BZ4}}]\label{cor:sen2}
	Let $(X,K)$ be a continuous-time BRW and
	consider a sequence of continuous-time BRWs $\{(X_n, {K_n})\}_{n\in \N}$
	such that $\liminf_{n \to \infty} X_n = X$.
	Suppose that $k_{xy}(n)  \le k_{xy}$ for all $n\in\N$, $x,y\in X_n$
	and $k_{xy}(n) \to k_{xy}$ as $n \to \infty$ for all $x,y\in X$.
	Then $\lambda_s^{(X_n, {K_n})}(x_0)\ge \lambda_s^{(X,K)}(x_0)$ for every $n \in \mathbb{N}$ and $x_0 \in X_n$. Moreover, given $x_0 \in X$, 
	we have $x_0 \in X_n$, eventually as $n \to +\infty$ and
	$\lambda_s^{(X_n, {K_n})}(x_0)
	\to \lambda_s^{(X,K)}(x_0)$ as $n\to \infty$.
\end{Corollary}

Among all possible choices of the sequence
$\{(X_n,\mu_n)\}_{n \in \N}$ there is one which is
\textit{induced} by $(X,\mu)$ on the subsets $\{X_n\}_{n\in \N}$;
more precisely, one can take $\mu_n(g):=\sum_{f \in S_X:f|_{X_n}=g} \mu_x(f)$
for all $x \in X_n$ and $g \in S_{X_n}$.
Roughly speaking, this choice means that all reproductions outside $X_n$ are suppressed.
In this case it is simply
$m(n)_{xy}=m_{xy}$ for all $x,y \in X_n$.

\section{Strong local survival and local modifications}
\label{sec:survivalprob}

In this section, we examine how modifications on finite sets may 
influence the behaviour of a BRW. Since 
survival on a fixed set for a continuous-time BRW is equivalent to 
survival on the same set for its discrete-time counterpart, it is 
natural to use results from the discrete-time case to infer 
corresponding results for continuous-time BRWs.

We recall here \cite[Theorem 4.1]{cf:BBHZ}, which is a fundamental tool for comparing extinction probability vectors.
In the case of global survival (choose $B=X$), it gives equivalent conditions for strong local survival
in terms of extinction probabilities.
We denote by ${\mathbf{q}}_0({x},A)$ the probability that the process, which starts with one particle at $x$, has no progeny ever in the set $A$.
Theorem~\ref{th:strongconditioned} compares the extinction probabilities of two sets. Specifically, it asserts that a BRW has a positive probability of surviving in a set $B$ while becoming extinct in a set $A$ if and only if it has a positive probability of surviving in $B$ without ever visiting $A$. The ``if'' direction of this equivalence is straightforward, whereas the converse is not. This is a very powerful result, which has many applications.

\begin{Theorem}[{\cite[Theorem 4.1]{cf:BBHZ}}]\label{th:strongconditioned}
	For any BRW 
	and $A,B \subseteq X$, the following statements are equivalent:
	\begin{enumerate} 
		\item there exists $x \in {X}$ such that ${\mathbf{q}}(x,B) < {\mathbf{q}}({x},A)$;
		\item there exists $x \in {X}$ such that ${\mathbf{q}}({x},B\setminus A) < {\mathbf{q}}({x},A)$;
		\item there exists $x \in {X}$ such that ${\mathbf{q}}({x},B) < {\mathbf{q}}_0({x},A)$;
		\item there exists $x \in {X}\setminus A$ such that,
		starting from $x$
		there is a positive chance of survival in $B$ without ever visiting $A$;
		\item there exists $x \in {X}$ such that, starting from $x$ there is a positive chance of survival in $B$ and extinction in $A$;
		\item 
		\[
		\inf_{x \in {X}\colon {\mathbf{q}}(x,B)<1} \frac{1-{\mathbf{q}}(x,A)}{1-{\mathbf{q}}(x,B)}=0.
		\]
	\end{enumerate}
\end{Theorem}
From Theorem~\ref{th:strongconditioned}, which is stated for a single BRW, we derive Theorem~\ref{th:modifiedBRW} and
Corollary~\ref{cor:pureweak-nonstrong} which give us information about the behaviour
of a BRW after some modifications.
To this aim, consider two BRWs  $(X,\mu)$ and $(X,\nu)$
and denote by ${\mathbf{q}}^\mu$ and ${\mathbf{q}}^\nu$ their respective
extinction probability vectors. If $\mu_x=\nu_x$ for all $x \not \in A$, for 
some set $A\subseteq X$, we derive properties of ${\mathbf{q}}^\nu$
from analogous properties of ${\mathbf{q}}^\mu$.
Theorem~\ref{th:modifiedBRW} extends  \cite[theorem 4.2]{cf:BZ17}, which addresses the special case $B = X$.

\begin{Theorem}\label{th:modifiedBRW}
	Let $(X,\mu)$ and $(X,\nu)$ be two BRWs. Suppose that  $A \subseteq X$ is a nonempty set
	such that $\mu_x=\nu_x$ for all $x \not \in A$.
	Then we have that
	${\mathbf{q}}_0^\mu(x,A)={\mathbf{q}}_0^\nu(x,A)$ for all $x \in X$ and, for all $B \subseteq X$, 
	\[
	{\mathbf{q}}^\mu(A) \leq {\mathbf{q}}^\mu(B) \Longleftrightarrow
	{\mathbf{q}}^\nu(A) \leq {\mathbf{q}}^\nu(B).
	\]
	If, in addition,  $\mu_x=\nu_x$ for all $x \not \in B$, then
	\[
	{\mathbf{q}}^\mu(A) = {\mathbf{q}}^\mu(B) \Longleftrightarrow
	{\mathbf{q}}^\nu(A) = {\mathbf{q}}^\nu(B).
	\]
	%
	\end{Theorem}
	
	It is natural to compare the extinction probabilities of a set $A$ with those of another generic set $B$ for both processes. According to the previous theorem, one BRW has a positive probability of surviving in $B$ and going extinct in $A$ if and only if the same holds for the other. In particular, if the differences between the two BRWs are confined to $A \cap B$, then the extinction probabilities of $A$ and $B$ coincide for one process if and only if they coincide for the other. We emphasize that this does not imply that all four extinction probabilities are identical; for example, it may happen that ${\mathbf{q}}^\mu(A) = {\mathbf{q}}^\mu(B)=\mathbf{1}$ (almost sure extinction in $A$ and $B$), while ${\mathbf{q}}^\nu(A) = {\mathbf{q}}^\nu(B)<\mathbf{1}$ (the same positive probability of survival in both $A$ and $B$).
	
	\begin{proof}[Proof of Theorem~\ref{th:modifiedBRW}]
	We note that $(X,\mu)$ and $(X,\nu)$ have the same behaviour until they first hit $A$,
	hence ${\mathbf{q}}_0^\mu(x,A)={\mathbf{q}}_0^\nu(x,A)$ for all $x \not \in A$. If $x \in A$, then clearly
	${\mathbf{q}}_0^\mu(x,A)=0={\mathbf{q}}_0^\nu(x,A)$.
	\\
	Suppose now that ${\mathbf{q}}^\mu(A) \nleq {\mathbf{q}}^\mu(B)$; this is equivalent to $\mathbf{q}(x,A)>\mathbf{q}(x,B)$ for some $x \in X$. Hence, according to Theorem~\ref{th:strongconditioned},
	there exists  $x \in X\setminus A$ such that there is a positive probability of survival in $B$ starting from $x$ without ever visiting $A$.
	Since the two processes have the same behaviour until they first hit $A$, the same holds for $(X,\nu)$ and this implies that 
	${\mathbf{q}}^\nu(x,A)>{\mathbf{q}}^\nu(x, B)$; thus
	${\mathbf{q}}^\nu(A) \nleq  {\mathbf{q}}^\nu(B)$. The reversed implication follows by exchanging the role of $\mu$ and $\nu$.
	\\
	The final double implication, 
	$
	{\mathbf{q}}^\mu(A) = {\mathbf{q}}^\mu(B) \Longleftrightarrow
	{\mathbf{q}}^\nu(A) = {\mathbf{q}}^\nu(B),
	$
	follows from the previous one by exchanging the role of $A$ and $B$.
%
\end{proof}

In particular,
Note that, if $A_0:=\{x \in X \colon \mu_x \neq \nu_x\}$, then the partial order of the sets $\{\mathbf{q}^i(A) \colon A_0 \subseteq A \subseteq X\}$ ($i=\mu, \, \nu$) is preserved by the (well-defined) map $\mathbf{q}^\mu(A) \mapsto \mathbf{q}^\nu(A)$.

Theorem \ref{th:modifiedBRW}, applied to continuous-time BRWs, leads to Corollary~\ref{cor:pureweak-nonstrong}, which describes how 
local modifications of the BRW affect the critical parameters.
This corollary is a slight extension of the first part of \cite[Corollary 4.4]{cf:BZ17}; specifically, we added the final line to the statement.

\begin{Corollary}[{\cite[Corollary 4.4]{cf:BZ17}}]\label{cor:pureweak-nonstrong}
Let $(X,K)$ and $(X,K^*)$ two irreducible continuous-time BRWs such that $k_{xy}=k^*_{xy}$ for all $x \in X \setminus A$
where $A$ is a finite nonempty set. 
Then the following are equivalent:
\begin{enumerate}
\item $\lambda^*_w<\lambda_w$;
\item $\lambda^*_s<\lambda_w$;
\item $\lambda^*_w=\lambda^*_s<\lambda_w$.
\end{enumerate}
In particular, if $\min(\lambda_s , \lambda_s^*) \ge \max(\lambda_w, \lambda_w^*)$ (for instance, if $\lambda_s =\lambda_s^*$), then $\lambda_w=\lambda_w^*$.
\end{Corollary}

This corollary leads to the following significant conclusion: either the weak critical parameter $\lambda_s$ of a BRW attains its maximum value (under finite modifications), or the BRW does not exhibit a pure global survival phase.

\begin{proof}[Proof of Corollary~\ref{cor:pureweak-nonstrong}]
Observe that the discrete-time counterparts of these continuous-time BRWs satisfy the hypotheses
of Theorem~\ref{th:modifiedBRW}, namely, their offspring distribution are the same outside $A$.

Clearly 
$2. \Longrightarrow 1.$ and $3. \Longrightarrow 2.$.
We just need to prove that $1. \Longrightarrow 3.$; more precisely, we prove that 
$\lambda^*_w<\lambda_w \Longrightarrow \lambda^*_w=\lambda^*_s$. Take $\lambda \in (\lambda_w^*, \lambda_w)$;
the $\lambda$-$(X,K^*)$ BRW survives globally by definition, hence ${\mathbf{q}}^*(X) < \mathbf{1}$ 
(we denote by ${\mathbf{q}}^*$ the extinction probabilities of $(X,K^*)$).
On the other hand, since $\lambda<\lambda_w$ the $\lambda$-$(X,K)$ BRW dies out a.s.~in every set starting from every point, whence  $\mathbf{1}={\mathbf{q}}(X)=\mathbf{q}(A)$. According
to Theorem~\ref{th:modifiedBRW}, since ${\mathbf{q}}(X)=\mathbf{q}(A)$ then ${\mathbf{q}}^*(X)=\mathbf{q}^*(A)$, which implies
$\mathbf{q}^*(A)< \mathbf{1}$. If the $\lambda$-$(X,K^*)$ BRW survives locally in
the finite set $A$ it means that it survives locally at a vertex $x \in A$ ($\Longleftrightarrow$ at every vertex, since
the process is irreducible). This implies $\lambda \ge \lambda_s^*$. We just proved that for every $\lambda \in (\lambda_w^*, \lambda_w)$ we have BRW survives globally, thus $\lambda_s^* \le\lambda_w^*$, whence $\lambda_s^* = \lambda_w^*$. 
\end{proof}

According to Corollary~\ref{cor:pureweak-nonstrong}, when 
$(X,K)$ and $(X,K^*)$ are two irreducible continuous-time BRWs such that $k_{xy}=k^*_{xy}$ for all $x \in X \setminus A$
where $A$ is a finite nonempty set, then the only alternatives are:
\begin{enumerate}
\item $\lambda_w=\lambda_w^* \le \min(\lambda_s,\lambda_s^*)$,
\item $\lambda_w=\lambda_s < \lambda_w^* \le \lambda_s^*$,
\item $\lambda_w^*=\lambda_s^* < \lambda_w \le \lambda_s$.
\end{enumerate}
We are thus able to describe the behaviour of the modified BRW, when the original BRW has a pure global survival phase and the modified has no weak phase. The phase diagram can be seen in Figure \ref{fig:modifiedBRW}. Note that, while in general a BRW may survive or die out globally at its global critical parameter, here we know the behaviour at $\lambda=\lambda^*_w$ (the process dies out globally). Moreover, for $\lambda>\lambda_s$ there is local survival which is strong or non-strong depending on the behaviour of $(X,K)$.
These results are summarized in Corollary \ref{cor:pureweak-nonstrong2}.

\begin{Corollary}[{\cite[Corollary 4.4]{cf:BZ17}}]\label{cor:pureweak-nonstrong2}  
Let $(X,K)$ and $(X,K^*)$ two irreducible continuous-time BRWs such that $k_{xy}=k^*_{xy}$ for all $x \in X \setminus A$
where $A$ is a finite nonempty set.  Suppose that $\lambda^*_w<\lambda_w < \lambda_s$, and let $B\subseteq X$ be nonempty.
Then, for the BRW $(X,K^*)$
\begin{enumerate} 
\item if $\lambda \le \lambda_w^*$ there is a.s.~extinction in $B$;
\item if $\lambda \in (\lambda_w^*, \lambda_w)$ there is strong survival in $B$; 
\item if $\lambda =\lambda_w$ and $(X,K)$ dies out globally, then there is strong survival in $B$;
conversely, if $(X,K)$ survives globally and $B$ is finite; there is non-strong survival in $B$;
\item if $\lambda \in (\lambda_w, \lambda_s]$ and $B$ is finite, then
there is non-strong survival in $B$;
\item if $\lambda > \lambda_s$ and $B$ is finite, then survival is strong (resp.~non-strong) in $B$
if and only if 
the same holds for $(X,K)$.  
\end{enumerate}
\end{Corollary}

According to this corollary, if a BRW exhibits a pure global survival phase, then any local modification without such a phase necessarily displays a non-strong survival phase. The critical behaviour at $\lambda=\lambda_s$ for a general BRW is well understood: almost sure local extinction occurs (see \cite[Theorem 4.7]{cf:BZ2}). In contrast, at the global critical point $\lambda=\lambda_w$, the process may either survive (see \cite[Example 3]{cf:BZ2}) or become extinct almost surely (see, for instance, \cite[Theorem 4.8]{cf:BZ2}). However, if a BRW is obtained as a finite modification of another BRW, then at the global critical value $\lambda=\lambda_s$ extinction occurs almost surely, as stated in the first point of the previous corollary.
\begin{figure}[h!]
\begin{tikzpicture}[scale=1]
\draw[->] (0,-.8) -- (11,-.8);
\draw[-] (0,-.9) -- (0,-.7);
\draw[-] (5,-.9) -- (5,-0.7);
\draw[-] (9,-.9) -- (9,-.7);

\node[above] at (-1,0.2) {$(X,K^*)$ $\rightarrow$};
\node[below] at (-1,-1) {$(X,K)$ $\rightarrow$};
\draw [decorate,decoration={brace,amplitude=5pt,mirror,raise=4ex}]
(0,-.3) -- (5,-.3) node[midway,yshift=-3em]{global ext.};
\draw [decorate,decoration={brace,amplitude=5pt,mirror,raise=4ex}]
(5,-.3) -- (9,-.3) node[midway,yshift=-3em]{pure global phase};
\draw [decorate,decoration={brace,amplitude=5pt,mirror,raise=4ex}]
(9,-.3) -- (11,-.3) node[midway,yshift=-3em]{local surv.\phantom{g}};  
\draw[->] (0,0) -- (11,0);
\draw[-] (0,-.1) -- (0,.1);
\draw[-] (2,-.1) -- (2,.1);
\draw[-] (5,-.1) -- (5,.1);
\draw[-] (9,-.1) -- (9,.1);

\node[below] at (0,-0.15) {0};
\node[below] at (2,-0.1) {$\lambda_w^*=\lambda_s^*$};
\draw [decorate,decoration={brace,amplitude=5pt}]
(0,.2) -- (2,.2) node[midway,yshift=1em]{global ext.};

\node[below] at (5,-0.1) {$\lambda_w$};
\draw [decorate,decoration={brace,amplitude=5pt}]
(2,.2) -- (5,.2) node[midway,yshift=1em]{strong local surv.};

\node[below] at (9,-0.1) {$\lambda_s$};
\draw [decorate,decoration={brace,amplitude=5pt}]
(5,.2) -- (9,.2) node[midway,yshift=1em]{non-strong local surv.};

\draw [decorate,decoration={brace,amplitude=5pt}]
(9,.2) -- (11,.2) node[midway,yshift=1em]{local surv.};
\end{tikzpicture}
\caption{
Phase diagram of $(X, K)$ and its modification $(X, K^*)$, observed on a finite set $B \neq \emptyset$. The process $(X, K)$ has a pure global survival phase (i.e., $\lambda_w < \lambda_s$), while $(X, K^*)$ has a lower global critical parameter and no pure global survival phase: $\lambda_w^* = \lambda_s^*$.
When $\lambda < \lambda_w$, the modified process $(X, K^*)$ exhibits strong local survival. For $\lambda_w < \lambda \le \lambda_s$, it exhibits non-strong local survival. 
For $\lambda > \lambda_s$, the two processes behave identically: one exhibits strong local survival if and only if the other does.
}
\label{fig:modifiedBRW}
\end{figure}
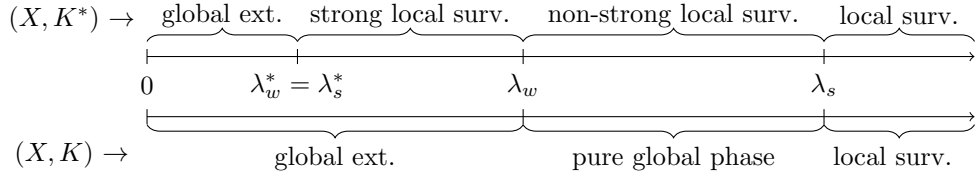
\begin{proof}[Proof of Corollary~\ref{cor:pureweak-nonstrong2}] 
We consider the following disjoint intervals for $\lambda$.
\begin{enumerate} 
\item Suppose that $\lambda < \lambda_w^*$: by definition, there is global extinction, hence also local extinction. 
From Corollary \ref{cor:pureweak-nonstrong2}, we know that $\lambda_w^* =
\lambda_s^*$.
If $\lambda = \lambda_w^*$, by \cite[Theorem 4.7]{cf:BZ2}, for any 
finite set $C$
the $\lambda$-$(X,K^*)$ BRW dies out locally at $C$. This means that
$\mathbf{q}^*(C)= \mathbf{1}$. We also know that $\mathbf{q}(B)={\mathbf{q}}(X)= \mathbf{1}$ for all $B\subseteq X$), since $\lambda < \lambda_w$.
Note that $\mathbf{q}(X)=\mathbf{q}(A)\le {\mathbf{q}}(B)$ and apply 
Theorem~\ref{th:modifiedBRW} to deduce
\[
{\mathbf{q}}^*(X)=\mathbf{q}^*(A)=\mathbf{1} \le \mathbf{q}^*(B).
\]
Thus, we have $\mathbf{q}^*(B)=\mathbf{1}$.

\item $\lambda \in (\lambda_w^*, \lambda_w)$. By definition, since $\lambda_w^*=\lambda_s^*$,
there is global and local survival for the $\lambda$-$(X,K^*)$ BRW. This implies that 
${\mathbf{q}}^*(X) \le \mathbf{q}^*(B) < \mathbf{1}$ for every (finite or infinite) set $B$. On the other hand,
there is global and local extinction for the $\lambda$-$(X,K)$ BRW, which 
implies that
${\mathbf{q}}(X)= \mathbf{q}(B)= \mathbf{1}$. Again, according to Theorem~\ref{th:modifiedBRW},
${\mathbf{q}}^*(X) = \mathbf{q}^*(B)$ for every  set $B$.
\item Clearly, since $\lambda =\lambda_w \le \lambda_s$, we have $\mathbf{q}(B)= \mathbf{1}$ for all finite subsets $B$.
Hence
\[
\lambda-(X,K) \textrm{ survives globally} \Longleftrightarrow {\mathbf{q}}(X)< \mathbf{q}(A),
\]
that is, according to Theorem~\ref{th:modifiedBRW}, if and only if ${\mathbf{q}}^*(X) < \mathbf{q}^*(A)$.
This, again, implies ${\mathbf{q}}^*(X) < \mathbf{q}^*(B)$ for every nonempty finite subset $B$.
If, on the other hand, $\lambda-(X,K)$  dies out globally, then ${\mathbf{q}}^*(X) = \mathbf{q}^*(A)$ and
${\mathbf{q}}^*(X) = \mathbf{q}^*(B)$ for every nonempty subset $B$.
\item $\lambda \in (\lambda_w, \lambda_s]$ (we suppose that the interval is nonempty, otherwise there is nothing to prove). 
Here we have ${\mathbf{q}}(X)< \mathbf{1}=\mathbf{q}(B)$ for every finite subset $B$.
Theorem~\ref{th:modifiedBRW} yields ${\mathbf{q}}^*(X) < \mathbf{q}^*(A)< \mathbf{1}$  and, by irreducibility,
${\mathbf{q}}^*(X) < \mathbf{q}^*(B)< \mathbf{1}$, for every finite nonempty subset $B$.
\item $\lambda > \lambda_s$. Now, $\mathbf{q}(B)< \mathbf{1}$ and $\mathbf{q}^*(B)< \mathbf{1}$
for every nonempty $B \subset X$.
Again, by Theorem~\ref{th:modifiedBRW}, we have
\[
{\mathbf{q}}(A)= {\mathbf{q}}(X) \Longleftrightarrow
{\mathbf{q}}^*(A)={\mathbf{q}}^*(X).
\]
If $B$ is finite, then Theorem~\ref{th:modifiedBRW} yields the conclusion.
\end{enumerate}
\end{proof}
In particular, if we know that $(X,K)$ exhibits non-strong local survival in a set $A$, for a fixed value
of $\lambda$, we can deduce the behaviour of a BRW, which coincides with $(X,K)$
outside $A$.

\begin{Corollary}\label{cor:pureweak-nonstrong3}
Let $(X,K)$ and $(X,K^*)$ two irreducible continuous-time BRWs such that $k_{xy}=k^*_{xy}$ for all $x \in X \setminus A$
where $A$ is a nonempty (not necessarily finite) set.
If there exists $\lambda_0$ such that $(X,K)$ survives non-strong locally in $A$, then 
\begin{enumerate}
\item $(X, K^*)$ survives globally when $\lambda=\lambda_0$;
\item if $A$ is finite and $\lambda^*_s \ge \lambda_0$, then there is pure global survival for $(X, K^*)$, if $\lambda=\lambda_0$.
\end{enumerate}
\end{Corollary}
\begin{proof}
If $\lambda=\lambda_0$,
$\mathbf{q}(A) >\mathbf{q}(X)$, then, according to Theorem~\ref{th:modifiedBRW},
$\mathbf{q}^*(A) >\mathbf{q}^*(X)$; this implies (1) that $\mathbf{q}^*(X) < \mathbf{1}$ (there is global survival) and (2) if $\lambda^*_s \ge \lambda_0$ and $A$ is finite, then $\mathbf{1}=\mathbf{q}(\{x\})=\mathbf{q}(A) >\mathbf{q}(X)$, where $x$ is a generic point (there is a pure global survival phase).
\end{proof}

A simple class of irreducible BRWs where $\lambda_w=\lambda_s$ is described by the following result. Here we consider the Alexandroff one-point compactification of $X$ endowed by the discrete topology: the neighbourhoods of $\infty$ are $\{X \setminus A \colon A \subset X, \textrm{ finite}\}$. The $\limsup_{x \to \infty}$ is defined accordingly.

\begin{Proposition}\label{pro:lambdaequal}
%
%
Let $(X,K)$ be an irreducible BRW such that 
$\lambda_s \limsup_{x \to +\infty} \sum_y k_{xy}\le 1$; then $\lambda_s=\lambda_w$.
\end{Proposition}
%
%
%

\begin{proof}
By hypothesis, there exists $\alpha \ge 0$ such that (1) $\lambda_s \alpha \le 1$ and (2) for every $\varepsilon>0$ there exists a finite $A=A_\varepsilon \subseteq X$ such that $\sup_{x \in X \setminus A} \sum_y k_{xy} <\alpha+\varepsilon$.
Suppose, by contradiction, that $\lambda_w<\lambda_s$ and let $\lambda \in (\lambda_w, \lambda_s]$; for every $\varepsilon>0$, since there is a.s.~extinction in every finite set, by Theorem~\ref{th:strongconditioned} there is a positive probability of global survival starting from some $x \in X \setminus A_\varepsilon$ without visiting $A_\varepsilon$. Since the total number of particles generated outside $A_\varepsilon$ is dominated by a branching process with a reproduction rate $\lambda (\alpha+\varepsilon)$, we have $\lambda (\alpha+\varepsilon) >1$. Since this is true for every $\varepsilon >0$ whe have proven that $\lambda \in (\lambda_w, \lambda_s]$ implies $\lambda \alpha \ge 1$.
We have two cases:
\begin{enumerate}
\item if $\alpha=0$, then $\lambda \alpha <1$ therefore the set $(\lambda_w, \lambda_s]$ is empty, that is, $\lambda_s=\lambda_w$;
\item if $\alpha>0$, then $\lambda_s \alpha \le 1$ and $\lambda \in (\lambda_w, \lambda_s]$ imply $\lambda_w \alpha \ge 1$ which, in turn, implies $\lambda_w \ge \lambda_s$; thus $\lambda_w=\lambda_s$.
\end{enumerate}
\end{proof}

\section{Maximality of the pure global survival phase}
\label{sec:max}

Let us consider the set of irreducible continuous-time BRWs on $X$, namely $\mathcal{X}:=\{(X, K) \colon K=(k_{xy})_{x,y \in X} \text{ is irreducible}\}$.
Define the relation $\mathcal{R}$ on $\mathcal{X}$ by
$(X,K) \, \mathcal{R}\, (X,K^*)$ if and only if there exists a finite $A\subset X$ such that $k_{xy}=k^*_{xy}$ for all $x \in X \setminus A$ and all $y \in X$.
An equivalent definition is $(X,K) \, \mathcal{R}\, (X,K^*)$ if and only if the set $\Delta_{K, K^*}:=\{x \in X \colon \exists y \in X, \, k_{xy}\neq k^*_{xy}\}$ is finite.

Clearly $\Delta_{K,K}=\emptyset$, whence $\mathcal{R}$ is reflexive; moreover $\Delta_{K, K^*}=\Delta_{K^*, K}$ therefore $\mathcal{R}$ is symmetric. Finally, for all $K, K^*, \widetilde K$ we have $\Delta_{K, K^*} \subseteq \Delta_{K, \widetilde K} \cup \Delta_{\widetilde K, K^*}$ which implies that $\mathcal{R}$ is transitive. The relation $\mathcal{R}$ is an equivalence relation. To avoid a cumbersome notation, we denote by $[K]$ the equivalence class of $(X,K)$ and by $\mathcal{X}/_\mathcal{R}$ the quotient set.  Henceforth,
given an irreducible BRW $(X,K)$ we denote by $\lambda_w^K$ and $\lambda_s^K$ its global and local critical values respectively.

\begin{Proposition}\label{pro:maximality}
Let $(X,K) \in \mathcal{X}$ such that $\lambda_w^K < \lambda_s^K$. Then for all $K^* \in [K]$ we have $\lambda_w^{K^*} \le \lambda_w^K$. Moreover,  for all $K^* \in [K]$ such that $\lambda_w^{K^*} < \lambda_s^{K^*}$,  we have $\lambda_w^{K^*} = \lambda_w^K$.	
\end{Proposition}
\begin{proof}
By contradiction, if $\lambda_w^{K^*} > \lambda_w^K$, then, 
according to Corollary~\ref{cor:pureweak-nonstrong}, we have
$\lambda_s^K=\lambda_w^K$, which contradicts the hypothesis.
From the previous part, it easily follows that if $\lambda_w^K < \lambda_s^K$ and $\lambda_w^{K^*} < \lambda_s^{K^*}$, then $\lambda_w^{K^*} = \lambda_w^K$.
\end{proof}

We note that when $X$ is finite, there is only one equivalence class: in this case, no BRW admits a pure global survival phase (see Remark~\ref{rem:finiteX}), and the supremum of $\lambda_w$ is clearly infinite. Hence, the case of real interest arises when $X$ is infinite.

As an application of Proposition \ref{pro:maximality}, we discuss the behaviour of the critical parameters of the BRW on the $d$-dimensional regular tree $\mathbb{T}_d$, 
when we add the possibility of reproduction from the origin of the tree to itself. 

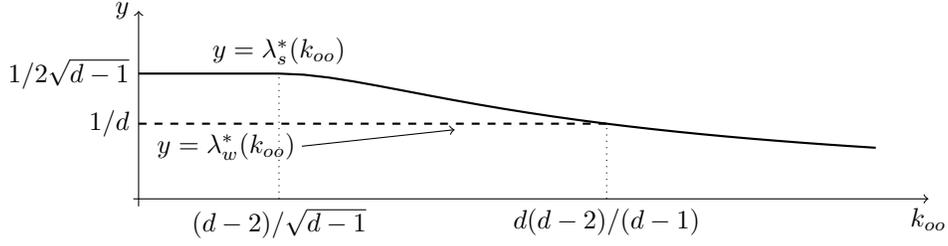
\begin{figure}[ht!]
\begin{tikzpicture}[xscale=0.7,yscale=10]
\draw[->] (-0.1,0) -- (15,0);
\draw[->] (0,-0.01) -- (0,0.25);
\node[left] at (0,0.25) {$y$};
\node[below] at (15,0) {$k_{oo}$};
\node[left] at (0,0.166667) {$1/2\sqrt{d-1}$};
\node[left] at (0,0.1) {$1/d$};

\node[below] at (2.666667,0) {$(d-2)/\sqrt{d-1}$};
\node[below] at (8.888889,0) {$d(d-2)/(d-1)$};

\node[above] at (2.666667,0.166667) {$y=\lambda_s^*(k_{oo})$};
\node[below] at (1.666667,0.1) {$y=\lambda_w^*(k_{oo})$};
\draw[->] (3.1,0.07) -- (6, 0.092);

\draw[-,thick] (0,0.166667) -- (2.666667,0.166667);
\draw[-, dashed, thick] (0,0.1) -- (8.888889,0.1);

\draw[dotted] (2.66667,0) -- (2.666667,0.166667);
\draw[dotted] (8.888889,0) -- (8.888889,0.1);

\draw[domain=2.66667:14, thick]  plot(\x,{(4*\x+5*sqrt(4+\x^2))/(100+9*\x^2)});

\end{tikzpicture}
\caption{
The critical parameters $\lambda_w^*$ and $\lambda_s^*$ for the process $(X, K)$ on the tree $\mathbb{T}_d$.
The rates are defined as $k_{xy} = 1$ for all $x \sim y$, and $k_{oo} \in [0, +\infty)$. The functions $y = \lambda_w^*(k_{oo})$ (dashed) and $y = \lambda_s^*(k_{oo})$ (solid) coincide for $k_{oo} \ge d(d - 2)/(d - 1)$.
Moreover, $y = \lambda_w^*(k_{oo}) = \lambda_w^*(0)$ for $k_{oo} \le d(d - 2)/(d - 1)$, and $y = \lambda_s^*(k_{oo}) = \lambda_s^*(0)$ for $k_{oo} \le (d - 2)/\sqrt{d - 1}$.
}
\label{fig:modifiedtree}
\end{figure}

\begin{Example}\label{exmp:modifiedtree}
Consider the BRW on $X=\mathbb{T}_d$, where $k_{xy}=1$ for all neighbouring couples $(x,y)$. 
For this BRW,
$\lambda_s=1/2\sqrt{d-1}$ and $\lambda_w=1/d$ and the process exhibits a pure global survival phase. 
We modify this BRW by adding a loop at the origin $o$, and choosing $k_{oo}>0$. We denote by $\lambda_w^*$ and $\lambda_s^*$ the global and local critical parameters, respectively, of this modified process. It is well known that for small values of $k_{oo}$, $\lambda_w^*<\lambda_s^*$, while for $k_{oo}$ 
sufficiently large, the critical parameters coincide  (see for instance \cite[Example 4.2]{cf:BZ14-SLS}). 
We are now able to prove that for $k_{oo}\in[0,(d-2)/\sqrt{d-1}]$, the two parameters coincide with the parameters of the original process ($k_{oo}=0$); for $k_{oo}\in((d-2)/\sqrt{d-1},d(d-2)/(d-1))$ we have $\lambda_w^*=\lambda_w$, while $\lambda_s^*\in(\lambda_w,\lambda_s)$ and is monotonically decreasing.
When $k_{oo}\ge d(d-2)/(d-1))$ the two parameters coincide and tend to 0 as $k_{oo}$ goes to infinity. The plot of the two functions $\lambda_w^*=\lambda_w^*(k_{oo})$ and $\lambda_s^*=\lambda_s^*(k_{oo})$ can be seen in Figure~\ref{fig:modifiedtree} (note that the axis do not share the same scale).

In order to compute the function $\lambda_s^*=\lambda_s^*(k_{oo})$, recall that the local critical parameter equals $\max\{\lambda >0 \colon \Phi(o,o|\lambda) \le 1\}$ (see Section \ref{sec:basic} or \cite{cf:BZ2}), where
\[
\Phi(o,o|\lambda)
= d \frac{1-\sqrt{1-4(d-1)\lambda^2}}{2(d-1)}+
{k_{oo}}\lambda
\]
The computation of the generating function $\Phi$ in this case is simple: there is only one path from $o$ to $o$ of length 1 (the loop) whence $\varphi^{(1)}_{oo}=k_{oo}\lambda$; on the other hand, the loop is not involved in any 1st-return path of length $n\ge 1$ therefore the coefficient $\varphi^{(n)}_{oo}$ for $n \ge 1$ is the same as in the homogeneous tree (that is, when $k_{oo}=0$). Hence, $\Phi$ is the sum of $k_{oo}\lambda$ and the analogous generating function of the homogeneous tree.
It is an exercise to obtain that, for $k_{oo}>(d-2)/\sqrt{d-1}$,
\[
\lambda_s^*=\frac{(d-2)k_{oo}+d\sqrt{k_{oo}^2+4}}{2(d^2+(d-1)k_{oo}^2)}.
\]
Note that this provides an example in which there are several (more precisely, infinitely many) local modifications of a BRW, all with $\lambda_w^*<\lambda_s^*$, and by Proposition \ref{pro:maximality} they all share the same global critical parameter. Moreover, one of these modifications (the one with $k_{oo}=d(d-2)/(d-1)$), has no pure global survival phase but still the same global critical parameter.
\end{Example}

Example~\ref{exmp:modifiedtree} was first presented in \cite[Example 4.2]{cf:BZ14-SLS}, where the authors proved that if the loop rate $k_{oo}$ is sufficiently large, then $\lambda^*_w=\lambda^*_s$. The proof in \cite{cf:BZ14-SLS} relies on a detailed analysis of the behaviour of a BRW along the branches of a homogeneous tree. In our case, we can instead employ results such as Corollary~\ref{cor:pureweak-nonstrong} and Proposition~\ref{pro:maximality}, which apply in more general settings than the one considered in \cite[Example 4.2]{cf:BZ14-SLS}. This not only simplifies the proof compared to \cite[Example 4.2]{cf:BZ14-SLS}, but also enables us to provide a complete characterization of the critical parameters of $(X, K^*)$ as functions of $k_{oo}$.


In the following sections we address several natural questions concerning these equivalence classes. We briefly summarize the main results (see also Table~\ref{tb:table2}). 
A given class may or may not contain BRWs with a pure global survival phase (see Section~\ref{subsec:q1}); however, every class necessarily contains uncountably many BRWs without a pure global survival phase (see Section~\ref{subsec:q2}). 
If a class admits at least one BRW with a pure global survival phase, then there are uncountably many BRWs with a pure global survival phase and the supremum of $\lambda_w$ is finite and attained by all such BRWs (see Section~\ref{subsec:q8} and Proposition~\ref{pro:maximality}). In this case, the same class may also contain BRWs without a pure global survival phase but with the same maximal value of $\lambda_w$ (see Section~\ref{subsec:q4}). 
If a class contains only BRWs without a pure global survival phase, then the supremum of $\lambda_w$ may be infinite (see Section~\ref{subsec:q6}) or finite; in the latter case, the supremum may or may not be attained (see Sections~\ref{subsec:q3} and \ref{subsec:q5}, respectively). 
Finally, we show that there exist classes in which $\lambda_w$ is constant across all members (see Section~\ref{subsec:q7}).

\begin{table}[h]
\setlength\extrarowheight{3pt}
\begin{tabular}{ |c|c| }
\hline 
\rowcolor{lightgray}
\multicolumn{2}{|c|}{\textbf{Classes without BRW with a pure global survival phase:
		$\lambda_w^*=\lambda_s^*$ for all $K^*\in [K]$}}   \\ \hline
\rowcolor{lightgray!50}
\textbf{Result} &  \textbf{Section}\\ 
\hline
such classes exist & \ref{subsec:q1}\\ 
\hline
it is possible that $\sup_{K^*\in [K]}\lambda_w^*<+\infty$ is attained & \ref{subsec:q3} \\ 
\hline
it is possible that $\sup_{K^*\in [K]}\lambda_w^*<+\infty$ and it is not attained & \ref{subsec:q5}  \\ 
\hline
it is possible that $\sup_{K^*\in [K]}\lambda_w^*=+\infty$ & \ref{subsec:q6}  \\ 
\hline
it is possible that $\lambda_w^*$ is constant on $K^*\in [K]$& \ref{subsec:q7}  \\ 
\hline
\hline
\rowcolor{lightgray}
\multicolumn{2}{|c|}{\textbf{Classes with a BRW with a pure global survival phase: $\lambda_w<\lambda_s$}}   \\ \hline
\rowcolor{lightgray!50}
\textbf{Result} &  \textbf{Section}\\ 
\hline
there is at least one $K^*\in[K]$ such that $\lambda_w^*=\lambda_s^*$& \ref{subsec:q2}  \\ 
\hline
it is possible that  
$\lambda_w=\lambda_w^*=\lambda_s^*$ for some $K^*\in[K]$
& \ref{subsec:q4}  \\ 
\hline
there exist uncountably many $K^*\in[K]$ such that $\lambda_w=\lambda_w^*<\lambda_s^*$& \ref{subsec:q8}  \\ 
\hline
\end{tabular}
\captionof{table}{
Summary of Section~\ref{sec:max}. 
We analyse the behaviour of the critical parameters $\lambda_w^*$ and $\lambda_s^*$ of $(X, K^*)$ as $K^*$ varies in $[K]$, and compare them to the critical parameters $\lambda_w$ and $\lambda_s$ of $(X, K)$.
}\label{tb:table2}
\end{table}

\subsection{Does every equivalence class have a BRW with a pure global survival phase?}\label{subsec:q1}
The answer is negative. Let us start with a heuristic reasoning before moving on to a more formal proof. Take $(\mathbb{Z}, K)$ where $k_{xy}:=1/2$ if $|x-y| = 1$ and $0$ otherwise. Consider a modification restricted to an interval $I$. Conditioned to global survival either an infinite number of particles visits $I$ or an infinite number of particles visits $\mathbb{Z} \setminus I$ but in this case, by a Borel-Cantelli argument, the progeny of an infinite number of particles visits $I$. 
So there is no pure global survival phase.

This can be proven rigorously in general as follows.
We recall that $(X,K)$ is \textit{quasi-transitive} if  there exists a finite $X_0\subset X$ such that for any $x\in X$ there is
a bijective map $\gamma_x:X\to X$ satisfying $\gamma_x^{-1}(x)
\in X_0$ and $k_{yz}=k_{\gamma_x y\,\gamma_x z}$ for all
$y,z$. For example, if the rates are translation invariant, then $(X,K)$ is quasi-transitive (actually, it is transitive).

\begin{Lemma}\label{lem:nopureglobalphase}
Let $(X,K)$ be a quasi-transitive, irreducible BRW on an infinite set $X$.
Suppose that for all $x \in X$ the set $\{y \in X\colon (x, y) \in E_\mu\}$ is finite.
If $(X,K)$ has no pure global survival phase, then every BRW in the same class has no pure global survival phase.
\end{Lemma}
\begin{proof}
Here are the main steps of the proofs. Given a subset $A \subseteq X$, we denote by $\lambda_w(A)$ and $\lambda_s(A)$ the critical parameters of the BRW $(A,K_A)$ where $K_A$ is the restriction of $K$ to $A \times A$ (this is the BRW restricted to $A$). Clearly $\lambda_w(A) \ge \lambda_w$ and $\lambda_s(A) \ge \lambda_s$. M
Moreover, according to Theorem \ref{th:spatial},
if $A_n \subseteq A_{n+1}$ and $\bigcup_n A_n=X$, then $\lambda_s(A_n) \downarrow \lambda_s$.
Due to the quasi-transitivity of the BRW, there is $x_0 \in X$ such that $\{\gamma(x_0)\colon \gamma \in \mathrm{AUT}(X,K)\}$ is infinite; clearly, for every automorphism $\gamma$ we have $\lambda_w(\gamma(A_n))=\lambda_w(A_n)$, $\lambda_s(\gamma(A_n))=\lambda_s(A_n)$. Let us choose $A_n$ finite for every $n$; observe that, for every fixed finite $I \subset X$, since $X$ is infinite and the BRW is quasi-transitive, there exists an automorphism $\gamma_n$ such that $\gamma_n(A_n) \cap I=\emptyset$. 
Indeed, let $r:=\max\{d(x_0, y)\colon y \in A_n\}$ where $d(x,y):=\min\{n \colon x \stackrel{n}{\to} y\}$ (which is not, in general, a distance because it might not be symmetric, but it satisfies the triangular inequality). Let the automorphism $\gamma_n$ be such that $d(\gamma_n(x_0),I)> r$, whose existence is guaranteed by the fact that $X$ is infinite and the set $\{y \in X\colon (x, y) \in E_\mu\}$ is finite for all $x \in X$. Thus, $\gamma_n(A_n) \cap I=\emptyset$.
\\
Let $K' \in [K]$ be such that the matrices are the same outside $I$. Clearly, if $\lambda_s' \le \lambda_s$ or $\lambda_w' \le \lambda_w$, then $\lambda_s'=\lambda_w'$ by Corollary~\ref{cor:pureweak-nonstrong}. Let us suppose that
$\lambda_w=\lambda_s \le \lambda_w' \le \lambda_s'$. Clearly, $\lambda_s(A_n) \ge \lambda_s'$ since $\gamma_n(A_n) \cap I=\emptyset$ and $\lambda_s(\gamma_n(A_n))=\lambda_s(A_n)$. Thus  
\[
\lambda_w=\lambda_s \le \lambda_w' \le \lambda_s' \le \lambda_s(A_n) \to \lambda_s, \qquad \textrm{as } n \to \infty
\]
therefore $\lambda_w=\lambda_s = \lambda_w' = \lambda_s'$.
\end{proof}


\subsection{Does every equivalence class have a BRW with no pure global survival phase?}\label{subsec:q2}

The answer is affirmative. Consider a BRW $(X,K)$ and denote by $\lambda_w$ its global critical parameter. Let $o \in X$ be fixed and define $K'$ by
\[
k'_{xy}:=k_{xy}+\alpha \delta_{xo}\delta_{yo}.
\]
Roughly speaking, we modify the rate from $o$ to itself by adding $\alpha$ to $k_{oo}$. Clearly $\lambda_s^\prime \le 1/\alpha$ therefore by taking $\alpha > 1/\lambda_w$ we have $\lambda^\prime_s<\lambda_w$ thus,
according to Corollary~\ref{cor:pureweak-nonstrong}, $\lambda^\prime_w=\lambda^\prime_s$.

\subsection{Is there a class without BRWs with a pure global survival phase where the supremum of $\lambda_w$ is attained?}\label{subsec:q3} 

The answer to this question is affirmative: it is possible to construct a BRW with no pure global survival phase, such that all BRWs in its equivalence class have no pure global survival phase, and where the supremum of $\lambda_w$ is attained. Consider the class of the BRW $(\mathbb{Z}, K)$, where $k_{xy}=1/2$ if $|x-y|=1$ and $0$ otherwise. 
It is straightforward to show that $1=\lambda^K_w=\lambda^K_s$: the first equality is immediate, while the second follows, for instance, from \cite[Theorem 4.3]{cf:BCZ}. Observe that if we restrict the same BRW to $\mathbb{N}$, the equalities $\lambda_w(\mathbb{N})=\lambda_s(\mathbb{N})=1$ still hold; this follows from \cite[Theorem 4.3]{cf:BCZ}, Theorem~\ref{th:spatial}, and from the fact that the critical parameters of a restricted BRW cannot be smaller than those of the original BRW. By Lemma~\ref{lem:nopureglobalphase}, there are no BRWs with a pure global survival phase in this class.

We now prove that the supremum of $\lambda_w$ in the equivalence class of $(\mathbb{Z}, K)$ is $1$.  
Let $K^* \in [K]$ and suppose that the modification is contained in $[-n,+n]$. Then, by a coupling argument, we obtain $\lambda_w^{K^*} = \lambda_s^{K^*}\le \lambda_w(\mathbb{N})=1$, since $(\mathbb{Z},K^*)$ dominates its restriction to $[n+1,+\infty)$, which is equivalent to the restriction of $(\mathbb{Z},K)$ to $\mathbb{N}$.  
Finally, to show that $\lambda_w^{K^*}$ is not constant within the class, note that by adding a loop at $0$ with sufficiently large rate $k_{00}$, one can construct a BRW in this class whose critical values satisfy $\lambda_w=\lambda_s<1$.

\subsection{If a class has no BRWs with a pure global survival phase and the supremum of $\lambda_w$ is finite, is it always attained?}\label{subsec:q5}

The answer is negative. Consider the class of $(\mathbb{Z}, K)$ where 
$k_{xy}=k_{|x|}/2$ if $|x-y|=1$ and
$k_n \downarrow \bar k$ ad $n \to \infty$. In this case the supremum of $\lambda_w$ equals $1/\bar k$, indeed take $K^*$ as
$k^*_{xy} :=\min(k_{xy}, (\bar k+\varepsilon)/2)$ where $\varepsilon>0$; since $\varepsilon$ is a generic positive constant, we have $\lambda_w^* \ge 1/(\bar k+\varepsilon)$ (since the restriction of a BRW to a subset gives a BRW with larger critical parameters). Then the supremum is at least $1/\bar k$. The reversed inequality will be proven in the specific example we are going to construct. The sequence $\{k_n\}_{n \in \mathbb{N}}$ will consist of long stretches of finite consecutive constant rates $\bar k+1/n$. More specifically, according to Theorem \ref{th:spatial},
for any $n \in \mathbb{N}$ there is $l_n \in \mathbb{N}$ such that the BRW on $\mathbb{Z}$ with rates $k_{xy}:=(\bar k+1/n)/2$ if $|x-y| =1$ and $0$ otherwise (where $\lambda_s =\lambda_w=1/(\bar k+1/n)$), restricted to a set $\{i+1,\ldots,i+l_n\}$ has a local critical parameter $\lambda_s^n$ close to $1/(\bar k+1/n)$, say $\lambda_s^n< 1/\bar k$. This does not depend on $i$.
Define $\{i_n\}$ recursively as $i_{n+1}:=i_n+l_{n+1}$, where $i_0:=0$, and $k_i:=\bar k +1/(n+1)$ for all $i=i_n, i_n+1, \ldots, i_{n+1}-1$. Consider $K^* \in [K]$; clearly there exists $n \in \mathbb{N}$ such that $k^*_i:=\bar k +1/(n+1)$ for all $i=i_n, i_n+1, \ldots, i_{n+1}-1$, therefore $\lambda_w^*\le \lambda_s^* \le \lambda_s^n<1/\bar k$. This shows that the supremum of $\lambda_w$ in class $[K]$ is $1/\bar k$ and is never attained.

\subsection{Is the supremum of $\lambda_w$ necessarily finite in a class?}\label{subsec:q6}

The answer is negative: clearly here we are dealing with classes with no BRWs with pure global survival phase. Take $(\mathbb{Z},K)$ where $k_{nm}=1/2(n^2+1)$ if $|n-m|=1$ and $0$ otherwise. This BRW satisfies the hypotheses of Proposition~\ref{pro:lambdaequal}, where $\limsup$ equals $0$. Consider the following modification in $[-n,+n]$, $k^*_{ij}=1/2(n^2+1)$ for all $i \in [-n,n]$ and $|i-j|=1$. This process is dominated by a branching process with an expected number of children equal to $\lambda/(n^2+1)$ therefore $\lambda^{K^*}_w \ge n^2+1$. Thus, the supremum of $\lambda_w$ in this class is infinite.

\subsection{Can $\lambda_w$ be constant in a class?}\label{subsec:q7}
The answer is affirmative. Clearly, according to Section~\ref{subsec:q2} and Corollary~\ref{cor:pureweak-nonstrong}, the class has no BRWs with a pure global survival phase. Take $(\mathbb{N}, K)$ where $k_{xy}=k_{|x|}/2$ if $|x-y|=1$ and
$k_n \uparrow \bar k$ ad $n \to \infty$. By spatial approximation (\cite[Theorem 4.3]{cf:BCZ}, Theorem \ref{th:spatial}
and domination, $1/\bar k \le \lambda_w^K \le \lambda_s^K \le 1/\bar k$. Since by spatial approximation $\lambda_s^{K^*}=1/k$ for any modification, then, according to Corollary~\ref{cor:pureweak-nonstrong}, $\lambda_w^{K^*}=\lambda_w^K$ (since $\lambda_s^{K^*} \ge \lambda_w^K$).

\subsection{If there is a BRW with a pure global survival phase in a class, can a BRW with no pure global survival phase have the same (maximum) $\lambda_w$?}\label{subsec:q4}
The answer is affirmative. We rely on Example~\ref{exmp:modifiedtree}.
It suffices to note that the BRW on $\mathbb T_d$ where $k_{xy}=1$ if $x$ and $y$ are neighbours and $k_{xy}=0$ otherwise, has a pure global survival phase,
while the modified BRW which adds a loop at the origin $o$, of rate $k_{oo}=\frac{d(d-2)}{d-1}$, has the same global critical parameter $1/d$, but no global
survival phase.

\subsection{Is there a class containing a positive, finite number of BRWs with a pure global survival phase?}\label{subsec:q8}
The answer is negative. We already know that there exist classes without BRWs exhibiting a pure global survival phase. Suppose that $(X,K)$ has a pure global survival phase.  
Select a vertex $o \in X$ and set the rate $k^*_{oo} = k_{oo} + \varepsilon$. Then, by the same arguments used in Example~\ref{exmp:modifiedtree}, we have 
$
\Phi^*(o,o|\lambda) = \Phi(o,o|\lambda) + \varepsilon \lambda.
$ 
By equation~\eqref{eq:lambdas1}, $\lambda_w < \lambda_s$ if and only if there exists some $\lambda_0 > \lambda_w$ such that $\Phi(o,o|\lambda_0) < 1$. Hence, for every sufficiently small $\varepsilon > 0$, we have $\Phi^*(o,o|\lambda_0) < 1$.  
All these uncountably many BRWs with a pure global survival phase share the same global critical parameter $\lambda_w$, according to Proposition~\ref{pro:maximality}.

\section{Conclusions}

In this paper, we investigate the effect of modifying the rates of a continuous-time BRW in a finite subset of the underlying space.
We say that two BRWs are equivalent if one can be obtained from the other through a local modification.

Consider a BRW $(X, K)$ that exhibits a pure global survival phase; that is, $\lambda^K_w < \lambda^K_s$.  
We prove that all other BRWs in its equivalence class that also have a pure global survival phase satisfy $\lambda_w^* = \lambda_w^K$:  
the local critical parameter may change, but the global one does not.
Moreover, $\lambda_w^K$ is the maximal global survival critical parameter among all BRWs in the equivalence class.  
Indeed, let $(X, K^*)$ be a local modification of $(X, K)$, and let $\lambda_w^*$ denote its global survival critical parameter. Then,  
$\lambda_w^* \le \lambda_w^K$.
If $\lambda_w^* < \lambda_w^K$, then $(X, K^*)$ does not exhibit a pure global survival phase.  
In this case, $\lambda_w^K$ becomes the threshold for $(X, K^*)$ between strong local survival in finite sets and non-strong local survival.

We analyse the existence of a pure global survival phase.
Given any BRW, its equivalence class always contains uncountably many processes that do not exhibit a pure global survival phase.
However, the converse does not hold: there exist equivalence classes that contain no BRWs with a pure global survival phase.

We provide several examples illustrating the possible behaviour of $\lambda_w$ within an equivalence class.
The global survival critical parameter may remain constant across the class, or it may be unbounded. Even when it is bounded, its supremum may not be attained.

In this paper, we consider only modifications within finite sets.
Allowing modifications on infinite sets would, in principle, lead to a situation where all processes become equivalent.
To avoid this triviality, one could restrict attention to modifications on connected infinite sets whose complement is also infinite.
However, studying such cases would require the development of new tools.

Among the possible generalizations of this work, a natural direction is to study equivalence classes of general BRWs, without assuming irreducibility.
Caution is needed in the general case, as some equivalence classes may contain both irreducible and reducible processes.
Recall that, in general, the critical values may depend on the starting vertex.
A key step toward this goal is to remove the irreducibility assumption from Corollaries~\ref{cor:pureweak-nonstrong} and~\ref{cor:pureweak-nonstrong2}.
We believe this is feasible, since these corollaries rely on Theorems~\ref{th:strongconditioned} and~\ref{th:modifiedBRW}, which do not require irreducibility.

Addressing other natural questions will likely require different techniques.
It is already known that isolating a superspreader can increase the critical parameters and lead an epidemic to extinction.
In general, understanding which locations have the greatest influence on the process can be far from obvious.
In the context of global survival, some works study how generations move toward infinity by identifying the so-called ends of the process and analysing its limiting behaviour at the boundary (see \cite{cf:Cand1, cf:Cand2, cf:KW}).
We conjecture that the ends most likely to be reached are those whose neighbourhoods are the most influential.

Another active area of research is the study of BRWs in random environment
(see \cite{cf:Gantert2, cf:MP00, cf:MP03, cf:Muller1}), where the rates are not deterministic but are instead chosen independently according to a fixed distribution.
In this setting, a local modification can be interpreted as a perturbation of the stochastic homogeneity of the rates: for instance, one may fix the rates within a finite set $A$, while allowing them to remain i.i.d. outside $A$.
This corresponds to exerting control over the rates in $A$, while preserving randomness elsewhere.
If the random model is such that the critical parameters do not depend on the specific realization of the environment, it is natural to ask whether this invariance persists after the modification, and how the critical parameters behave in the modified system.

\section*{Acknoledgements}

The authors thank the three anonymous referees for useful comments on the previous version of the paper.

\end{document}